\documentclass[12pt]{article}
\usepackage{latexsym, amssymb, amsmath, amscd, amsfonts, epsfig, graphicx, colordvi}

\newtheorem{thm}{Theorem}[section]

\def\pf{\noindent{\it Proof.} }
\def\qed{\nopagebreak\hfill{\rule{4pt}{7pt}}
\medbreak}

\def\qed{\nopagebreak\hfill{\rule{4pt}{7pt}}
\medbreak}

\newenvironment{kst}
{\setlength{\leftmargini}{2.4\parindent}
\begin{itemize}
\setlength{\itemsep}{-0.5mm}} {\end{itemize}}

\begin{document}

\begin{center}

{\Large \bf   Jacobi's Identity and Synchronized  Partitions}\\[15pt]
\end{center}

\begin{center}
{ William Y. C. Chen}$^{1}$ \quad  and  \quad {Kathy Q. Ji}$^{2}$

   Center for Combinatorics, LPMC\\
   Nankai University, Tianjin 300071, P.R. China

   \vskip 1mm

   Email: $^1$chen@nankai.edu.cn, $^2$ji@nankai.edu.cn
\end{center}

\vskip 6mm \noindent {\bf Abstract.} We obtain a finite form of
Jacobi's identity and present a combinatorial proof based on the structure
of synchronized partitions.



\section{Introduction}

We adopt the common notation on partitions  and $q$-series as used
in \cite{and76, gas90}. The {\it $q$-shifted factorial} $(x; q)_n$
is defined by $(x; q)_0 = 1$ and for $n\geq 1$,
\[
(x; q)_n = (1 - x)(1 - qx) \cdots (1- q^{n-1}x).\] The {\it
$q$-binomial coefficient}, or the {\it Gauss coefficient}, is
given by
\begin{equation}\label{bin1}
{n\brack k} =\frac{(q;q)_n}{(q;q)_k(q;q)_{n-k}}, \quad \mbox{ for
}  0 \leq k \leq n.
\end{equation}
Note that the parameter $q$ is often omitted in the notation of
the Gauss coefficients.

This paper is concerned with Jacobi's  identity
\cite [p.257, Eq.(5)]{jac29}  (see also
\cite[Theorem 357]{wri65})
\begin{equation}\label{jac}
(q;q)^3_{\infty}=\sum_{k = 0}^{\infty}(-1)^k (2k+1)q^{k+1 \choose
2}.
\end{equation}
Note that the identity (\ref{jac}) can be deduced from
Jacobi's  triple product identity \cite[p.15]{gas90}
\begin{equation}\label{jac-tr}
(z;q)_{\infty}(q/z;q)_{\infty}(q;q)_{\infty}=
\sum_{k=-\infty}^{\infty}(-1)^kq^{k\choose 2}z^k.
\end{equation}
 Rewriting (\ref{jac-tr})
in the following form:
\begin{equation} \label{jac-tr2}
(z;q)_{\infty}(z^{-1}q,;q)_{\infty}(q;q)_{\infty}=\sum_{k=0}
^{\infty}(-1)^k(1-z^{2k+1})z^{-k}q^{k+1\choose2},
\end{equation}
then one obtains (\ref{jac}) from (\ref{jac-tr2}) by dividing both
sides  by $(1-z)$ and taking the limit $z\rightarrow1$. A
combinatorial proof of (\ref{jac}) has been found by Joichi and
Stanton \cite{joi89}.

Jacobi's identity  has many  applications.
Ramanujan \cite{ram19,ram21}  proved the partition congruences
modulo 5 and 7 by using  this  identity. Andrews
\cite{and84,and96}, Ewell \cite{ewe82},   and  Hirschhorn and
Sellers \cite{hir96} have used this identity in their studies of
congruence relations on partition functions. Jacobi's identity
(\ref{jac}) also plays a role in the study of representing an
integer as sum of squares, see, for example,
 Hirschhorn \cite{hir99}.

The purpose of this paper is to derive the finite form of Jacobi's
  identity and give a combinatorial proof.  We first give a finite form of Jacobi's
  identity \eqref{jac} by using MacMahon's finite form of
Jacobi's triple product identity. Then we
 give the definitions of synchronized  partitions and
 rooted synchronized  partitions and present two simple
 involutions on synchronized  partitions which
 imply a combinatorial proof of the finite form of Jacobi's identity.

\section{A Finite Form of Jacobi's  Identity}

We obtain the following finite form of
Jacobi's identity.

\begin{thm}\label{thm}   For $m, n\geq 0$, we have
\begin{equation}\label{fin-jac}
(q;q)_m(q;q)_n=\sum_{k=-n-1}^{m}(-1)^k k q^{k+1 \choose 2} {m+n+1
\brack n+k+1}.
\end{equation}
\end{thm}

\pf We begin with MacMahon's finite form of Jacobi's triple
product identity \cite[Vol.~II, \S323]{ma15} (see also
\cite{chu93,han99})
\begin{equation}\label{fin-tri}
(zq;q)_m(z^{-1};q)_n=\sum_{k=-n}^{m}(-1)^kq^{k+1 \choose 2}z^k{m+n
\brack n+k}.
\end{equation}
Substituting $n$ with $n+1$ in the above identity, we get
\begin{equation}\label{eqn3}
(zq;q)_m(z^{-1};q)_{n+1}=\sum_{k=-n-1}^{m}(-1)^kq^{k+1 \choose
2}z^k{m+n+1 \brack n+1+k}.
\end{equation}
Setting
$$f(z)=(zq;q)_{m}(z^{-1}q;q)_n,$$
then (\ref{eqn3}) becomes
$$(1-z^{-1})f(z)=\sum_{k=-n-1}^{m}(-1)^kq^{k+1 \choose
2}z^k{m+n+1 \brack n+1+k}.$$
Differentiating both sides respect to $z$, we get
\begin{align*}
z^{-2}f(z)+(1-z^{-1})f'(z)= \sum_{k=-n-1}^{m}(-1)^kkq^{k+1 \choose
2}z^{k-1}{m+n+1 \brack n+1+k}.
\end{align*}
Setting $z=1,$  one obtains (\ref{fin-tri}). \qed

Setting $n\rightarrow\infty$ and $m\rightarrow\infty$ in
\eqref{fin-jac}, by Tannery's theorem (see\cite[p.292]{tan04}), we
get
\[(q;q)^3_{\infty}=\sum_{k=-\infty}^{\infty}(-1)^k k
q^{k+1 \choose 2},
\]
which is equivalent to  Jacobi's   identity (\ref{jac}).

Replacing $m$ by $n$ in \eqref{fin-jac}, we obtain
 the following
identity.

\begin{thm} For $n\geq 0$, we have
\begin{equation}\label{fin-jac2}
(q;q)_n^2=\sum_{k=0}^n(-1)^k(2k+1) q^{k+1 \choose 2}  {2n+1 \brack
n+k+1}.
\end{equation}
\end{thm}

\section{Synchronized Partitions}

In this section, we give a combinatorial proof of the finite form
of Jacobi's identity (\ref{fin-jac}) by introducing the structures of
synchronized partitions and rooted synchronized partitions. Let us
recall some common terminology on partitions.  {\it A partition}
$\lambda$ of a positive integer $n$ is a finite weakly decreasing
sequence of positive integers
$\lambda_1,\,\lambda_2,\ldots,\,\lambda_r$ such that
$\sum_{i=1}^r\lambda_i=n$, denoted by
$\lambda=(\lambda_1,\,\lambda_2,\ldots,\,\lambda_r)$, where the
$\lambda_i$'s are called the parts of $\lambda$; the sum of parts
is called the {\it weight} of  $\lambda$,  denoted by $|\lambda|$;
the number of parts of $\lambda$ is called the {\it length} of
$\lambda$,  denoted by $l(\lambda).$

A pair of partitions $(\alpha, \beta)$  of the same length
is called a generalized Frobenius partition, see Andrews
\cite{and84}, Corteel and Lovejoy \cite{cor02}. In a more general
setting, a pair of partitions $(\alpha, \beta)$ that are not
necessarily of the same length is also called a generalized
Frobenius partition, see  Yee \cite{yee03, yee04}.

We now give the definitions of  {\it synchronized partitions} and
{\it rooted  synchronized partitions}. Assume that
$\alpha=(\alpha_1, \alpha_2, \ldots, \alpha_r)$ is a partition
with distinct parts and $\beta=(\beta_1, \beta_2, \ldots,
\beta_s)$ is also a partition with distinct parts under the
assumption  that the last part $\beta_s$ may be zero. Then a
synchronized partition is a representation of $(\alpha, \beta)$ as
a two-row array such that some $*$ symbols may be added at the end
of $\alpha$ or $\beta$ so that they are of the same length
depending on which is of smaller length. We may denote a synchronized partition
with underlying partitions $\alpha$ and $\beta$ by $S(\alpha, \beta)$, or simply
$(\alpha, \beta)$ is no confusion arises. The difference $r-s$ is
called the {\it discrepancy}  of the synchronized partition. A
synchronized partition with a positive discrepancy $k$ can be
represented as follows:
\[S (\alpha,\,\beta)=\left(\begin{array}{ccccccc}
\alpha_1&\alpha_2&\cdots&\alpha_s&\alpha_{s+1}&\cdots&\alpha_{s+k}\\
\beta_1&\beta_2&\cdots&\beta_s&*&*&*
\end{array}\right)\]
and  a synchronized  partitions  with a negative discrepancy $-k$
($k> 0$) can be represented as follows:
\[S (\alpha,\,\beta)=\left(\begin{array}{ccccccc}
\alpha_1&\alpha_2&\cdots&\alpha_r&*&*&*\\
\beta_1&\beta_2&\cdots&\beta_r&\beta_{r+1}&\cdots&\beta_{r+k}
\end{array}\right).\]
A synchronized partition with zero discrepancy can be simply
represented as a two-row array without any star added. A {\it
rooted synchronized partition} is defined as a synchronized
partition with a distinguished star symbol, which we denote by
$\bar{*}$. Clearly, a rooted
synchronized partition has an underlying synchronized partition
with nonzero discrepancy.

For example, there are five rooted synchronized partitions  of $2$
 :
$$\left(\begin{array}{c}
2\\
\bar{*}
\end{array}\right)\left(\begin{array}{cc}
1&\bar{*}\\
1&0
\end{array}\right)\left(\begin{array}{cc}
\bar{*}\\
2
\end{array}\right)\left(\begin{array}{cc}
\bar{*}&*\\
2&0
\end{array}\right)\left(\begin{array}{cc}
*&\bar{*}\\
2&0
\end{array}\right)$$

Let $\mathcal{S}_{m,n}$ denote the set of synchronized partitions
$S(\alpha, \beta)$ such that $\alpha_1\leq m$ and $\beta_1\leq n$,
and let   $\mathcal{R}_{m,n}$ be the set of rooted synchronized
partitions $S(\alpha, \beta)$ such that $\alpha_1\leq m$ and
$\beta_1\leq n$. Note that $\mathcal{S}_{m,n}$ are generated by
the set of pairs of partitions $(\alpha, \beta)$ under the same
condition by adding some stars to a row if it is of smaller length
than the other row.
A rooted synchronized  partition $S(\alpha,\beta)$
 is called
{\it degenerate}  if
\[S=\left(\begin{array}{cccccc}
\alpha_1&\cdots&\alpha_s&\alpha_{s+1}&\cdots&\alpha_{r}\\
\beta_1&\cdots&\beta_s&\bar{*}&\cdots&*
\end{array}\right)\]
or
\[S=\left(\begin{array}{ccccccc}
\alpha_1&\cdots&\alpha_r&*&\cdots&*&\bar{*}\\
\beta_1&\cdots&\beta_r&\beta_{r+1}&\cdots&\beta_{s}&0
\end{array}\right)\]
where $m\geq \alpha_1>\alpha_2>\ldots>\alpha_r \geq 1$ and $n\geq
\beta_1>\beta_2>\ldots>\beta_s\geq 1;$ otherwise $S(\alpha,
\beta)$ is called {\it non-degenerate}.

 It is easy to see that the generating function of
synchronized  partitions in $\mathcal{S}_{m,n}$ equals
\begin{equation}
\sum_{S(\alpha, \beta) \in \mathcal{S}_{m,n}}   q^{|\alpha|+|\beta|} =
 (-q;q)_m(-1;q)_{n+1}, \end{equation}
 and the generating function of  synchronized partitions
 in $\mathcal{S}_{m,n}$ without the zero part
equals
\begin{equation}
(-q;q)_m(-q;q)_n.
\end{equation}
On the other hand,  the generating function of
 synchronized
partitions
 in $\mathcal{S}_{m,n}$ with  a nonnegative discrepancy $k$  equals
  \begin{equation} \label{q1}
  q^{k+1 \choose 2}{m+n+1 \brack n+k+1},
  \end{equation}
and the generating function of   synchronized
partitions
 in $\mathcal{S}_{m,n}$ with a negative discrepancy  $-k$  equals
 \begin{equation} \label{q2}
 q^{-k+1\choose 2} {m+n+1 \brack n-k+1}= q^{k \choose 2}{m+n+1 \brack n-k+1}.
 \end{equation}
From (\ref{q1}) and (\ref{q2}) it follows that
the generating function for  rooted synchronized
 partitions in $\mathcal{R}_{m,n}$ equals
\[\sum_{k=0}^mkq^{k+1 \choose 2}{m+n+1 \brack n+k+1}+
\sum_{k=1}^{n+1}kq^{k \choose 2}{m+n+1 \brack n-k+1}.\]

 Let us define the sign of a rooted synchronized partition
$S(\alpha, \beta)$ as $(-1)^{\delta (S)}$, where $\delta(S)$ is
the number of stars in $S(\alpha, \beta)$ under the assumption
that a star with the bar in the top row is not counted. The sign
of a synchronized partition equals $(-1)^k$, where $k$ is the
discrepancy. We now give a sign reversing involution on the set of
non-degenerate rooted synchronized partitions.

\begin{thm}\label{main1} There is a sign reversing involution $\tau$ on the set of
non-degenerate rooted synchronized partitions of $p$ in
$\mathcal{R}_{m,n}$.
\end{thm}

\pf For a non-degenerate rooted synchronized  partition
$S(\alpha,\beta)\in \mathcal{R}_{m,n},$  we proceed to construct a
non-degenerate rooted synchronized partition $S(\alpha', \beta')$.
We consider the following two cases.
\begin{kst}
\item[Case 1:] The partition $\beta$ has a zero part.
\begin{kst}
\item  If $l(\alpha)> l(\beta)$, then replace the zero part by a
star $*$.

\item If $l(\alpha)< l(\beta)$, then delete the whole column of
the zero part.
\end{kst}
\item[Case 2:] The partition $\beta$ has no zero part.
\begin{kst}
\item If $l(\alpha)> l(\beta)$, then replace the first  `$*$' on
the bottom row  by a zero part.
 \item If $l(\alpha)< l(\beta)$, then add
a zero part along with a star on the top as a column.
\end{kst}
 \end{kst}
The above bijection can be illustrated as follows:
 \vskip -0.7cm
 \begin{align*}
 \begin{array}{ccc}
\left(\begin{array}{cccccc}
\alpha_1&\cdots&\alpha_{s}&a_{s+1}&\cdots&\alpha_{r}\\
\beta_1&\cdots&0&*&\bar{*}&*
\end{array}\right)& \stackrel{l(\alpha)>l(\beta)}{\text{\hbox to 40pt {\leftarrowfill}}}\hskip-1.5cm{\text{\hbox to 40pt {\rightarrowfill}}} &\left(\begin{array}{cccccc}
\alpha_1&\cdots&\alpha_s&\alpha_{s+1}&\cdots&\alpha_{r}\\
\beta_1&\cdots&*&*&\bar{*}&*
\end{array}\right),\\
 \end{array}
 \end{align*}
  \begin{align*}
 \begin{array}{ccc}
\left(\begin{array}{ccccccc}
\alpha_1&\cdots&\alpha_r&*&\bar{*}&*&*\\
\beta_1&\cdots&\beta_r&\beta_{r+1}&\cdots&\beta_{s-1}&0
\end{array}\right)&\stackrel{l(\alpha)<l(\beta)}{\text{\hbox to 40pt {\leftarrowfill}}}
\hskip-1.5cm{\text{\hbox to 40pt {\rightarrowfill}}}&\left(\begin{array}{ccccccc}
\alpha_1&\cdots&\alpha_r&*&\bar{*}&*\\
\beta_1&\cdots&\beta_r&\beta_{r+1}&\cdots&\beta_{s-1}
\end{array}\right) .\\
 \end{array}
 \end{align*}
It is easy to check that the above construction is a sign
reversing involution. \qed

Next, we give a bijection  between  degenerate rooted synchronized
partitions and   synchronized partitions without the zero part.

\begin{thm}\label{main2}
There is a sign preserving bijection  between the set of
degenerate rooted synchronized partitions of $p$ in
$\mathcal{R}_{m,n}$ and the set of  synchronized
partitions of $p$ in $\mathcal{S}_{m,n}$ that do not contain the zero part.
\end{thm}

\pf  For a degenerated rooted synchronized partition
$S(\alpha,\beta)$ in $\mathcal{R}_{m,n}$, we can construct a
 synchronized
 partition $S(\alpha',\beta')$ in $\mathcal{S}_{m,n}$ that do not contain the zero
 part.
\begin{kst}
\item[Case 1] If $l(\alpha)>l(\beta)$, then delete the bar to the
first `$*$' on the bottom row.
 \begin{align*}
 \begin{array}{ccc}
 \left(\begin{array}{cccccc}
\alpha_1&\cdots&\alpha_s&a_{s+1}&\cdots&\alpha_{r}\\
\beta_1&\cdots&\beta_s&\bar{*}&\cdots&*
\end{array}\right)& \longleftrightarrow &
\left(\begin{array}{cccccc}
\alpha_1&\cdots&\alpha_s&\alpha_{s+1}&\cdots&\alpha_{r}\\
\beta_1&\cdots&\beta_s&*&\cdots&*
\end{array}\right)
\end{array}
 \end{align*}

 \item[Case 2] If $l(\alpha)<l(\beta)$, then delete
a zero part on the bottom row along with a barred star on the top
row.
 \begin{align*}
 \begin{array}{ccc}
\left(\begin{array}{ccccccc}
\alpha_1&\cdots&\alpha_r&*&\cdots&*&\bar{*}\\
\beta_1&\cdots&\beta_r&\beta_{r+1}&\cdots&\beta_{s}&0
\end{array}\right)& \longleftrightarrow &
\left(\begin{array}{ccccccc}
\alpha_1&\cdots&\alpha_r&*&\cdots&*\\
\beta_1&\cdots&\beta_r&\beta_{r+1}&\cdots&\beta_{s}
\end{array}\right)
 \end{array}
 \end{align*}

\end{kst}
Clearly, the procedure is reversible and it preserves the signs.
\qed

We are now ready to give a combinatorial interpretation of finite
form  of Jacobi's   identity \eqref{fin-jac}. It is easy to see
that $(q;q)_n(q;q)_m$
is the generating function of  signed synchronized  partitions
$\mathcal{S}_{m,n}$ without the zero part. Note that
\begin{align*}
\sum_{k=-n-1}^{m}(-1)^k k q^{k+1 \choose 2} {m+n+1 \brack n+k+1}
&=\sum_{k=0}^{m}(-1)^k k q^{k+1 \choose 2} {m+n+1 \brack
n+k+1}\\
&+\sum_{k=1}^{n+1}(-1)^{k-1} k q^{k \choose 2} {m+n+1
\brack n-k+1}
\end{align*}
is the generating function of signed rooted synchronized
partitions in $\mathcal{R}_{m,n}.$ Combining  Theorem \ref{main1}
and Theorem
 \ref{main2}, we are led to a combinatorial interpretation of the
finite form \eqref{fin-jac}.

 \vspace{6mm}

\noindent{\bf Acknowledgments.}  This work was supported by the
973 Project on Mathematical Mechanization, the National Science
Foundation, the Ministry of Education, and the Ministry of Science
and Technology of China.

\end{document}